\documentclass[a4paper,draft,intlimits,reqno,10pt,oneside]{amsart}
\usepackage{amssymb,enumerate,eucal,verbatim}

\usepackage{amsmath}
\usepackage{amssymb}
\usepackage{amsthm}
\usepackage{latexsym}

\newtheorem{lemma}{Lemma}[section]
\newtheorem{theorem}[lemma]{Theorem}
\newtheorem{proposition}[lemma]{Proposition}
\newtheorem{corollary}[lemma]{Corollary}

\newtheorem{definition}[lemma]{Definition}
\newtheorem{fact}[lemma]{Fact}

\theoremstyle{definition}
\newtheorem{remark}[lemma]{Remark}

\def\({\bigl(}
\def\){\bigr)}

\def\vep{\varepsilon}

\def\N{I\!\!N}
\def\R{I\!\!R}

\def\fr#1#2{\dfrac{\lower 1 pt\hbox{\mathstrut}#1}{\mathstrut #2}}

\def\N{I\!\!N}

\def\R{I\!\!R}

\def\({\big(}
\def\){\big)}
\def\vep{\varepsilon}

\def\N{I\!\!N}
\def\R{I\!\!R}

\def\vs{\vskip 5 mm}
\def\vep{\varepsilon}
\def\EndDef{\end{definition}}
\def\Def{\begin{definition}}
\def\EndFact{\end{fact}}
\def\Fact{\begin{fact}}
\def\Pf{\begin{proof}}
\def\EndPf{\end{proof}}

\def\Lm{\begin{lemma}}
\def\EndLm{\end{lemma}}

\def\EndCor{\end{corollary}}
\def\Cor{\begin{corollary}}

\def\Thmn{\begin{theorem}}
\def\EndThm{\end{theorem}}
\def\Propn{\begin{proposition}}
\def\EndProp{\end{proposition}}

\newcommand\abs[1]{\left|#1\right|}

\newcommand\conv{\mathrm{conv}}

\def\N{I\!\!N}
\def\R{I\!\!R}

\def\vs{\vskip 5 mm}
\def\vep{\varepsilon}

\DeclareMathOperator{\dist}{dist}

\begin{document}

\title{Smooth approximations of norms in separable Banach spaces}

\author{Petr H\'ajek}
\address{Mathematical Institute\\Czech Academy of Science\\\v Zitn\'a 25\\115 67 Praha 1\\Czech Republic}
\email{hajek@math.cas.cz}

\author{Jarno Talponen}
\address{Aalto University, Institute of Mathematics, P.O. Box 11100, FI-00076 Aalto, Finland}
\email{talponen@iki.fi}

\date{\today}

\thanks{The first author was supported in part by Institutional Research Plan AV0Z10190503 and GA\v CR P201/11/0345.
This paper was prepared as the second author enjoyed the warm hospitality of the Czech Academy of Sciences in Spring
2011. The visit and research was supported in part by the V\"ais\"al\"a foundation.}

\begin{abstract}
Let $X$ be a separable real Banach space having a $k$-times continuously Fr\'{e}chet differentiable
(i.e. $C^k$-smooth) norm where $k\in \{1,\ldots,\infty\}$. We show that any equivalent norm on $X$ can
be approximated uniformly on bounded sets by $C^{k}$-smooth norms.
\end{abstract}

\maketitle

\section{Introduction}

The problem of approximation of continuous mappings by more regular functions, such as
polynomials or $C^k$-smooth functions, is one of the classical themes in analysis.
If the underlying spaces are infinite-dimensional the additional difficulty is in the
lack of compactness and measure as available tools. Yet, many important and rather
general results can be obtained even in this setting, as shown in the monograph
by Deville, Godefroy and Zizler \cite{dgz} and the references therein. In particular, if $X$ is
a separable real Banach space admitting a $C^k$-smooth bump function, $k\in\N\cup\{\infty\}$,
then every continuous mapping $f:X\to Y$, $Y$ a real Banach space, can be uniformly
approximated by means of $C^k$-smooth mappings \cite{bonfram}. The pioneering result in this area
is due to Kurzweil \cite{kurzweil} who proved a similar result in the real analytic setting
(provided $X$ has a separating polynomial), see also \cite{fry1}, \cite{ch}.

A natural version of the problem in the setting of norms is a question
whether every equivalent renorming of a (separable) real Banach space can
be uniformly approximated on bounded sets by means of $C^k$-smooth norms,
provided that $X$ admits an equivalent $C^k$-smooth norm.
Of course, from the above mentioned results we know that such approximations are
possible by means of $C^k$-smooth functions, but the additional requirement
of convexity cannot be achieved using these techniques.
This problem was posed explicitly in \cite{dgz} (for the case of a separable Hilbert space, see pp. 206-207). 
In \cite{defoha-separ} and \cite{defoha-polyhedral} the problem was solved in the affirmative for separable Banach spaces under some additional assumptions. The result was shown to hold if $X$ is either a polyhedral space (e.g. $c_0$), or
if it is a superreflexive space with a Schauder basis and the highest derivative $D^k\|\cdot\|$
is bounded (e.g. $\ell_2$)
. These papers left open not only the case of an abstract separable Banach space,
but also many concrete spaces such as $c_0\oplus\ell_2$ (which admits a $C^\infty$-smooth
renorming, but the $C^2$-smooth approximations of norms were open).

The main result of the present note is a complete solution of the problem
for all separable Banach spaces given in Theorem \ref{main-theo}.
The proof follows the main ideas in \cite{defoha-separ} with an essential new ingredient
contained in Lemma \ref{l: key}. The trick is to replace a given $C^k$-smooth norm on $X$
by a perturbation using a quadratic polynomial in order to get a new $C^k$-smooth renorming
which preserves $C^k$-smoothness with respect to quotients with a finite-dimensional kernel.
For the sake of convenience we will give a full proof of our theorem, providing also the arguments
for some known auxiliary results (e.g. Lemmas \ref{l: cantorset} and \ref{l: decomp}).
We also correct a minor mistake in the proof of the main results in \cite{defoha-separ},
which consisted of overlooking that the seminorms used in the proof are not smooth
at the points of their kernels. The exact place of this omission is indicated in our proof
of the main theorem.

\begin{comment}
In the final remarks we give some applications of our new technique to real algebraic geometry.
In particular, we show that the approximating Nash functions for a projection of a convex Nash function in the Efroymson Approximation Theorem can be chosen to be convex as well.
\end{comment}

\subsection{Preliminaries}

We denote by $X$, $Y$ and $E$ real Banach spaces. The closed unit ball and the unit sphere of $X$ are denoted by
$B_{X}$ and $S_{X}$, respectively. We will write $B(x,r)=x+rB_{X}$ for $x\in X$, $r>0$. For suitable background information on general Banach space theory and notations we refer to \cite{fhhmz}. However, next we will recall some basic concepts for the sake of convenience.

Recall that a mapping $f\colon X\to Y$ is \emph{Fr\'{e}chet differentiable} at point $x_{0}$ if there exists
a bounded linear mapping $T_{x_{0}}\colon X\to Y$ such that
\[\lim_{t\to 0^{+}}\sup_{h\in tS_{X}}\frac{\|f(x_{0}+h)-f(x_{0})-T_{x_{0}}(h)\|}{t}=0.\]
In such a case we denote by $\ Df(x_{0})[z]=T_{x_{0}}(z),$ the \emph{Fr\'{e}chet derivative}
of $f$ at $x_{0}$. We denote by $Df\colon X\to\mathcal{L}(X; Y)$
the Fr\'echet derivative as a function acting from the space $X$
with values in the space of bounded linear operators $\mathcal{L}(X;Y)$.
If $x\mapsto Df(x)$ is defined and continuous then $f$ is termed as being
\emph{continuously Fr\'{e}chet differentiable}, or in the class $C^{1}$.
In the case where $Df(x)$ is continuously Fr\'{e}chet differentiable with respect to $x$, we say that $f$ is (Fr\'{e}chet) $C^{2}$-smooth, and so on. If $f$ is $C^{k}$-smooth for all $k\in \N$, then it is $C^{\infty}$-smooth.
When we say that a norm is smooth we mean that it is that away from the origin.
As customary, we denote the class of real analytic functions by $C^{\omega}$. Here we will study such functions
only in the finite-dimensional setting. Given a map $f\colon X\oplus Y \to E$, we denote by $D_{2}f(x,y)$ the
Fr\'{e}chet differential taken with repect to the second coordinate $y$.

Given a subset $A\subset X$, a retraction $\rho\colon X\to A$ is a \emph{closest point mapping} if
$\|x-\rho(x)\|=\dist(x,A)$ holds for all $x\in X$. The closed linear span of $A$ is denoted by $[A]$.

Recall that a family $\{(x_{i},x_{i}^{\ast})\}_{i\in \N}$ of elements of $S_{X}\times X^{\ast}$ is a biorthogonal system if
$x_{i}^{\ast}(x_{j})$ equals to $1$ if $i=j$ and $0$ otherwise. It is called an \emph{$M$-basis} if
\[[x_{i}:\ i\in\N]=X\ \mathrm{and}\ \overline{[x_{i}^{\ast}:\ i\in\N]}^{\omega^{\ast}}=X^{\ast}.\]
If additionally $[x_{i}^{\ast}:\ i\in\N]=X^{\ast}$, then $\{(x_{i},x_{i}^{\ast})\}_{i\in \N}$ is a \emph{shrinking} $M$-basis.

Recall that a closed convex bounded (CCB) subset $C$ of a Banach space $X$ is a \emph{body} if it has a
non-empty interior. The \emph{polar} $C^{0}$ of $C$ is given by
\[C^{0}=\{f\in X^{\ast}:\ |f(x)|\leq 1,\ x\in C\}.\]
The following fact is well-known and can be found in \cite{dieudonne}.

\begin{theorem}[Implicit Function Theorem]
Let $X, Y , Z$ be Banach spaces, $U\subset X \oplus Y$ be an open set and $F\colon U \to Z$ be $C^1$-smooth.
Let $x_{0}, y_{0} \in U$ and $D_2 F(x_0 , y_0 )$ be a linear isomorphism of Y onto Z. Then
there is a neighbourhood $V$ of $x_0$ in $X$, and a unique mapping $u\colon V \to Y$,
such that $u(x_0 )= y_0$, and for all $x \in V$ it holds that
$(x, u(x)) \in U$ and $F(x, u(x))=F(x_0 , y_0 )$. Moreover, if $F$ is $C^k$-smooth,
$k\in \{1,2,\ldots,\infty,\omega\}$, then so is u.
\end{theorem}

\section{Results}
\subsection{Preparations}
\ \\
Recall some elementary facts from \cite{hp}.

\begin{definition}
We say that a function $f: \ell_\infty(\Gamma) \to \R$ is \emph{strongly lattice} if $f(x)\leq f(y)$
whenever $\abs{x(\gamma)} \leq \abs{y(\gamma)}$ for all $\gamma \in \Gamma$.
\end{definition}

\begin{definition}
Let $X$ be a vector space. A function $g:X \to \ell_\infty(\Gamma)$ is said to be
\emph{coordinate-wise convex} if, for each $\gamma \in \Gamma$, the function
$x \mapsto g_\gamma(x)$ is convex.
Similarly $g$ is \emph{coordinate-wise non-negative} if $g_\gamma(x)\geq 0$ for all
$\gamma \in \Gamma$.
\end{definition}

\begin{lemma}\label{l:SLconvexComposition}
Let $f:\ell_\infty(\Gamma) \to \R$ be convex and strongly lattice. Let
$g:X \to \ell_\infty(\Gamma)$ be coordinate-wise convex and coordinate-wise non-negative.
Then $f \circ g:X \to \R$ is convex.
\end{lemma}
\begin{proof}
Let $a,b \geq 0$ and $a+b=1$.
Since $g$ is coordinate-wise convex and non-negative, we have
\[
0\leq g_\gamma(ax+by)\leq a g_\gamma(x)+ b g_\gamma(y)
\]
for each $\gamma \in \Gamma$.
The strongly lattice property and the convexity of $f$ yield
\[
f(g(ax+by)) \leq f(a g(x)+ b g(y)) \leq a f(g(x))+ b f(g(y))
\]
so $f \circ g$ is convex.
\end{proof}

\begin{lemma}\label{l: key}
Suppose that $(X,\|\cdot\|)$ is a separable Banach space with a Fr\'echet $C^{k}$-smooth norm,
$k\in \{2,3,\ldots,\infty\}$, and $\{(x_{i},x_{i}^*)\}_{i\in \N}$ is an M-basis of $X$.
Then $\|\cdot\|$ can be approximated by norms $|||\cdot|||$ such that on $(X,|||\cdot|||)$ 
the quotient norms $|||\cdot |||_{X/[x_{1},\ldots,x_{n}]}$, $n\in\N$, are Fr\'echet $C^{k}$-smooth.
\end{lemma}
\begin{proof}
The required norms $|||\cdot|||$ are given by
\[|||x|||^{2}=\|x\|^{2}+\vep \sum_{i=1}^{\infty}2^{-i}x_{i}^{\ast}(x)^{2}, \quad \vep>0.\]
It is clear that these norms approximate the original one uniformly (on $\|\cdot\|$-bounded sets)
as $\vep$ tends to $0$. We observe immediately that the functional $|||\cdot|||^2$ is $C^{k}$-smooth
and therefore $|||\cdot|||$ is $C^k$-smooth away from the origin. From now on we will work with the space
$(X,|||\cdot|||)$.

Write $E=[x_{1},\ldots,x_{n}]$. Our strategy is to apply the Implicit Function Theorem to establish the  $C^k$-smoothness
of the (single-valued) closest point mapping $\rho\colon X \to E$.
Then the quotient norm $|||\cdot |||_{X/E}$ will be $C^k$-smooth away from the origin, since it can be written as
\[|||\hat{x} |||_{X/E}=|||x-\rho(x)|||,\]
which is necessarily $C^{k}$-smooth map as a composition of $C^k$-smooth mappings \mbox{$x-\rho(x)$} and
$|||\cdot|||$.

Fix $x\in X$ and let $d=\dist(x,E)$. By compactness the set of closest points
$\bigcap_{r>d}B(x,r)\cap E$ is non-empty. Suppose that it contains two points, say $a,b\in E$.
Then for $1\leq i\leq n$ we have that $x_{i}^{\ast}(x-\frac{a+b}{2})^{2}<\frac{x_{i}^{\ast}(x-a)^2 + x_{i}^{\ast}(x-b)^2}{2}$ unless $x_{i}^{\ast}(a-b)=0$. Thus $|||x-\frac{a+b}{2}|||<\frac{|||x-a|||+|||x-b|||}{2}$ if $a\neq b$. Hence the set of closest points is a singleton. In the sequel we shall denote by $\rho$ the closest point
mapping $X\to E$. Note that for $x\in X$ and $y\in E$ we have
\begin{eqnarray*}
\rho(x+y)&=&B(x+y,\dist(x+y,E))\cap E=y+(B(x,\dist(x+y,E))\cap E)\\
&=&y+(B(x,\dist(x,E))\cap E)=y+\rho(x).
\end{eqnarray*}

Next, we wish to verify that $\rho$ is $C^k$-smooth.
Let $L: X\oplus E\to X$ be a bounded linear mapping $L(x,y)=x-y$.

Let $G:X\oplus E\to (X\oplus E)^*$ be given by (using the chain rule)
\begin{equation*}
\begin{array}{lll}
G(x,y)&=&D(|||\cdot|||^2\circ L)(x,y)=D(|||\cdot|||^2)\circ L(x,y)\\
 &=&D|||x-y|||^2\in (X\oplus E)^*=X^*\oplus E^*.
\end{array}
\end{equation*}

Let $P: X^*\oplus E^*\to E^*$ be the canonical projection.
Note that since $\dim(E)=n$, there is a natural identification
$$
I: E^*\to E,\ I\left(\sum_{i=1}^n a_i x^*_i\right)=\sum_{i=1}^n a_i x_i.
$$
We let now $F: X\oplus E\to E$ be given by
$$
F(x,y)=I\circ P\circ G(x,y).
$$

Note that $F(x,y)=0$ if and only if $y=\rho(x)$. Indeed, $F(x,y)=0$ if and only if
$P\circ G(x,y)=0$ which means that $y\to |||x-y|||^2$
attains its minimum, due to the strict convexity of the last function.
Using the coordinates, if $y=\sum_{i=1}^n a_i x_i$ then

\[F(x,y)=\left(\frac\partial{\partial x_1}|||x+\sum_{i=1}^n a_i x_i|||^2, \dots,
\frac\partial{\partial x_n}|||x+\sum_{i=1}^n a_i x_i|||^2\right).\]

Recall that $|||\cdot|||^2$ is at least $C^2$-smooth. So, we may put

$$
D_2 F=DF|_E=\left(\frac{\partial^2}{\partial x_i\partial x_j}|||x+\sum_{i=1}^n a_i x_i|||^2\right)_{i,j=1,\dots, n}
$$

For a fixed $x$ the above map $D_2 F$ can be considered as a (symmetric) Hessian matrix
$M\in \R^{n \times n}$ of the restriction of $|||\cdot|||^2$ to the affine finite-dimensional
space $(x+E)$ in $X$. (See \cite[p.39]{borvein-vander} for the properties of the Hessian and convexity).

We claim that $M$ is invertible.  Indeed, invertibility follows easily from the the fact that $M$ is positive-definite. The functional $|||\cdot|||^{2}$ decomposes in a natural way into two parts, namely,
$\|x\|^{2}+\vep\sum_{n+1}^{\infty}2^{-m}x_{m}^{\ast}(x)^2$ and
$\vep\sum_{i=1}^{n}2^{-i}x_{i}^{\ast}(x)^2$.
The first part is a convex $C^{2}$-smooth function so that the corresponding Hessian matrix $M_{1}$ is
positive-semidefinite.
The Hessian matrix $M_{2}$ of the latter one is a strictly positive diagonal matrix, because of
the definition of the
biorthogonal functionals, thus positive-definite. By linearity, we obtain that $M=M_{1}+M_{2}$ is
positive-definite.
We conclude that $D_2 F(x,y)$ is an isomorphism $E\to E$ for any $x\in X, y\in E$.

Recall that $|||\cdot|||^2$ is a convex function so $F(x,y)=0$ exactly at the pairs of
points $(x,y)=(x,\rho(x))$. Fix $x_{0}\in X$. We apply the Implicit Function Theorem to find open
neighborhoods
$V\subset X$ of $x_{0}$ and $U\subset X\oplus E$ of $(x_{0},\rho(x_{0}))$, and a $C^{k}$-smooth mapping
$u\colon V\to E$ such that $F(x,u(x))=F(x_{0},\rho(x_{0}))=0$ for $(x,u(x))\in U$.
Necessarily $u(x)=\rho(x)$, as $\rho$ is single-valued. Since $x_{0}$ was arbitrary, we conclude that
$\rho$ is $C^{k}$-smooth on the whole of $X$.
\end{proof}

\begin{remark}
The smooth approximation of the quotient $X/E$ in the previous proof works for $E=\ell^{2}$ as well, essentially by the same argument.
\end{remark}

The following lemma is similar to Lemma 4.1 \cite{zp} and to the first part of
Theorem 1 \cite{ben}. We give a proof for the reader's convenience. Recall
that a Banach space $X$ with a separable dual admits a shrinking $M$-basis, see \cite[p.8]{bos}.

\vs

\begin{lemma}\label{l: cantorset}
Let $X$ be a Banach space with separable dual $X^*$ and let $\{(x_i, x_{i}^*)\}_{i\in\N}$ be a shrinking
$M$-basis of $X$. Let $W\subset X$ be CCB body such that  $\vec 0\in \mathrm{int}(W)$ and $0<\vep <\frac18$.
Then there exists a $\omega^*$-compact subset $F\subset W^0$ such that:
\begin{enumerate}
\item[(1)]{$\frac1{1+4\vep}W^0\subset \overline{\conv}^{\omega^{\ast}}(F)\subset \frac1{1+\vep}W^0$}
\item[(2)]{For each integer $i$ the set $x_i (F)$ is finite.}
\end{enumerate}
\end{lemma}

\begin{proof}
Let $d=\inf\{\|g\|:g\in \partial W^0 \}$ and $T_i=\{f(x_i):f\in
W^0\}\;\;\mathrm{for}\;\;i\in\N.$ Each set $T_i$ is bounded and thus there
exists a $\frac{\vep d}{6\cdot 2^{i}\|x_i^*\|}$-net $C_i$ in $T_i$. Put

$$ A=\left\{\sum\limits_{i=1}\limits^n a_i x_i^*\in\frac1{1+\vep}W^0 :\;n\in\N,\; a_i\in C_i \right\}, \;F=\overline{A}^{\omega^{\ast}}.$$

It is  obvious that $x_i (F)=x_i(A) \subset C_i$ where $i\in\N$. Thus condition 2.
is satisfied. Clearly $\overline{\conv}^{\omega^{\ast}}(F) \subset\frac1{1+\vep}W^0.$

Next, we will show that $\frac1{1+4\vep}W^0 \subset\overline{\conv}^{\omega^{\ast}}(F)$. Let us take $f\in \frac1{1+2\vep}W^0$ and recall that $\vep <\frac{1}{8}$.
Since $\frac1{1+4\vep}W^0$ is convex, by using the geometric Hahn-Banach theorem the task reduces to checking that there
exists $h\in F$ such that
\begin{equation}\label{eq: fheq}
\|f-h\|<d\frac{\vep}{3}<\frac{d\vep}{(1+\frac{1}{4})(1+\frac{1}{2})}<\frac{d2\vep}{(1+2\vep)(1+4\vep)}=d\left(\frac{1}{1+2\vep}-\frac{1}{1+4\vep}\right).
\end{equation}

Since
$\mathrm{span}\{x_i^*\}$ is dense in $X^*$, there exists
$$g=\sum\limits_{i=1}\limits^n b_ix_i^* \in\frac 1{1+2\vep}W^0, $$
such that $\|f-g\|<\vep\frac{d}6$. We have $b_i\in  T_i$ and there exists $a_i\in
C_i$  such that  $|b_i-a_i|< \frac{\vep d}{6\cdot 2^i \|x_i^*\|} , i\in\N.$
Hence
$$\|g-\sum\limits_{i=1}\limits^n a_i x_i^*\|<d\frac{\vep}{6},$$
and a straight verification shows that $h=\sum\limits_1\limits^n a_ix_i^* \in
\frac1{1+\vep}W^0$. Thus, $h\in F$ by the constructions of $h$ and the set $A\subset F$.
Now,
\[\|f-h\|\leq \|f-g\|+\|g-h\|<d\frac{\vep}{3},\]
so that \eqref{eq: fheq} holds, and the proof is complete.
\end{proof}

The following Lemma  is close to some results from [Zp] also. We will use
the notations from Lemma \ref{l: cantorset}. In addition put $M_n=[x_i]_1^{n \bot }, n\in\N.$

\begin{lemma}\label{l: decomp}
 For arbitrary $\vep >0$ there exists a sequence of
points $\{g_k\}$ in the set $F$, a sequence of integers $\{n_k\},n_k
\rightarrow \infty ,$ and a decreasing sequence $\{F_{\alpha}\}$ of
$\omega^*$-closed subsets of $F$ such that after reindexing it holds that
\begin{enumerate}
\item[(1)]{$\bigcup_{k\in\N} ((g_k +M_{n_k})\cap F_k)=F$.}
\item[(2)]{$\mathrm{diam}((g_k+M_{n_k})\cap F_k )< \vep$.}
\end{enumerate}
\end{lemma}

\begin{proof}
We will use the following well-known property of $\omega^*$-compacts in a separable dual space:
for every $\vep >0$ there exist a point $g\in F $ and $\omega^*$-neighborhood $G$ of $g$ such that $G\cap F \neq \emptyset $ and
$\mathrm{diam}(G\cap F)<\vep.$ Indeed, we apply the fact that $X$ is Asplund and $X^{\ast}$ has the RNP,
see \cite[Ch. 11.2]{fhhmz} for discussion.

Condition (2) in Lemma \ref{l: cantorset} suggests that in a sense $F$ resembles a Cantor set.
Because of the structure of the set $F$, the sets $(h+M_n)\cap F,h\in F, n\in\N$ form a base of
$\omega^*$- topology on $F$ and each such a set is both closed and open subset of $(F,\omega^*)$. Moreover, the family

$$\Im =\{h+M_n :h\in F, n\in N\}$$
contains countably many (different) sets and obviously each non-empty $\omega^*$- compact
subset $I\cap F$, $I\in \Im$, has the same structure as $F$. For each ordinal $\alpha$ we define sets $F_{\alpha}$ and
$(h_{\alpha}+M_{n(\alpha )})$ by transfinite induction as follows:

$$F_0=F,F_{\alpha +1}=F_{\alpha} \setminus (h_{\alpha}+M_{n(\alpha )}),$$
where $(h_{\alpha}+M_{n(\alpha)})$ is a member of the family $\Im$ such that
$(h_{\alpha}+M_{n(\alpha)})\cap F_{\alpha}\neq \emptyset $ and
$\mathrm{diam}((h_{\alpha}+M_{n(\alpha)})\cap F_{\alpha} )< \vep.$ If $\alpha$ is a limit
ordinal we put $F_{\alpha}=\cap_{\beta <\alpha }F_{\beta}.$ Since $F$ is separable and $\omega^*$-compact,
there exists a countable ordinal $\eta$ such that $F_{\eta}\neq \emptyset$ and
$F_{\eta+1}=\emptyset.$ It is clear that
$$\bigcup_{\alpha \leq \eta }\((h_{\alpha}+M_{n(\alpha)}\) \cap F_{\alpha})=F.$$
Let us reindex the countable family $\{h_{\alpha}+M_{n(\alpha)}\}_{\alpha\leq \eta }$
into $\{h_k+M_{n_k}\}_{k=1}^{\infty}.$ Since for each integer
$q$ there exist only finite many members $h+M_n$ of the family $\Im$  such
that $n\leq q$, it follows that $n_k \rightarrow \infty
\;\;\mathrm{for}\;\;k\rightarrow \infty .$
\end{proof}

\subsection{Fr\'{e}chet $C^k$-smooth approximation of norms}
\ \\
Next we give our main result.

\begin{theorem}\label{main-theo}
Let $(X,\|\cdot\|)$ be a separable Banach space. Let $k\in\N\cup\{+\infty\}$ and $\|\cdot\|$ be
$C^k$-smooth. Then every equivalent norm on $X$ can be approximated uniformly on bounded sets by $C^k$-smooth equivalent norms.
\end{theorem}

\begin{proof}

For $k=1$ the above result is known, see \cite[p.53]{dgz}.
 Since $X$ is Fr\'echet smooth the dual is separable
(see e.g. \cite[Thm. 8.6]{fhhmz}) and therefore $X$ admits a shrinking M-basis $\{(x_i, x_i^*)\}_{i\in\N}$.
By Lemma \ref{l: key} we may assume without loss of generality that the 
quotient norms $\|\cdot\|_{X/[x_{1},\ldots,x_{n}]}$, $n\in\N$, are Fr\'echet $C^{k}$-smooth.
Denote by $W=B_{(X,\|\cdot\|)}$,
$P_{n}\colon X \to [x_i]_{1}^{n}$ the projection given by $P_n(x)=\sum_{i=1}^{n}x_{i}^* (x)x_i$.
Put $Q_{n}(x)=x - P_{n}(x)$.
Using previous lemmas we obtain sequences $\{g_j\}_{j\in\N}\subset F$ and $\{M_{n_j}\}_{j\in\N}$ and put
\begin{equation}\label{eq: C}
C=\bigcup_{j\in\N}\(g_j+\vep\cdot B_{X^*}\cap M_{n_j}\)\supset F.
\end{equation}
By taking into account that $\vep$ is arbitrary, both above and in condition (1) of Lemma \ref{l: cantorset},
we observe that sets of the form $\overline{\conv}^{\omega^{\ast}}(C)$ are sufficient in approximating the polar $W^\circ$.
In fact, for technical reasons we will approximate $C$ with yet another set $D$ to be defined shortly.

Write $E_{j}=[x_i:\ 1\leq i\leq n_j]$ and $Y_{j}=[x_i:\ n_j +1\leq i]$ for $j\in \N$ and observe that
$X=E_j \oplus Y_j$, as $E_j$ is finite-dimensional. 
We will denote by $q_{n_{j}}$ the corresponding quotient mapping
$X\to X/E_j$ and we note that $(X/E_j)^{\ast}=(M_{n_j} ,\|\cdot\|^{\ast})$.
Observe that
\[\|q_{n_j}(x)\|_{X/E_j}=\|q_{n_j}(Q_{n_j}(x))\|_{X/E_j}\quad \mathrm{for}\ x\in X\]
by the definition of the quotient norm. Also, putting $q_{n_j}\circ Q_{n_j}$ defines an isomorphism $Y_j \to X/E_j$.
Unfortunately, the map $\|q_{n_j}(\cdot)\|_{X/E_j}$ is $C^k$-smooth only away from its kernel. 
This fact has been overlooked in \cite{defoha-separ} (where our $Q_n$ is denoted as $P_n$),
so the proof therein in not entirely correct. Therefore we need  to smoothen up this mapping.

Next we will give norms that are smooth away from the origin and  approximate the seminorms
$\|q_{n_j}(\cdot)\|_{X/E_j}$. Let $N$ be a norm on $\R^2$ satisfying the following conditions:
\begin{itemize}
\item[(a)]{$(\R^2,N)$ is $C^{\infty}$-smooth,}
\item[(b)]{$N(-a,b)=N(a,b)=N(b,a)$ for $a,b \in \R$,}
\item[(c)]{$(1,-1/2),(1,0),(1,1/2)\in S_{(R^2,N)}$.}
\end{itemize}
Indeed, this can be accomplished by taking a Minkowski functional of a suitable sum of two functions,
similarly as in \eqref{eq: F} to follow. Define norms $N_{j}\colon X = E_j \oplus Y_j \to [0,\infty)$ by
\[N_{j}(x) = N_{j}(P_{n_j}(x),Q_{n_j}(x))=N(a_{j}\|P_{n_j}(x)\| ,\|q_{n_{j}}(Q_{n_j}(x))\|_{X/E_{j}})\]
where $a_{j}\searrow0$ is a decreasing sequence of positive numbers. The main point here is that 
according to $(c)$ these norms coincide with the mapping $x\mapsto a_j\|P_{n_j}(x)\|$ around the kernel of 
$q_{n_j}$. Also, $N_j$ coincide with the mapping $x\mapsto \|q_{n_j}(Q_{n_j}(x))\|$ around the kernel of
$P_{n_j}$.
These facts together with the $C^k$-smoothness of all the mappings involved
in the definition are responsible for  $N_j$ being $C^k$-smooth away from the origin.
Note that $N_{j}(x)\geq \|q_{n_{j}}(x)\|_{X/E_{n_j}}$, so that
\[B_{X^*}\cap M_{n_j}\subset N_j^{\circ}\quad \mathrm{for}\ j\in\N.\]
Since $a_{j}\to 0$ as $j\to \infty$ we get
\begin{equation}\label{eq: (c)}
\mathrm{dist}_{\mathrm{H}}(N_j^{\circ},B_{X^*}\cap M_{n_j})\to 0\ \mathrm{as}\ j\to \infty.
\end{equation}
The above limit involves the symmetric Hausdorff distance.
Fix $\vep>0$ and put
\begin{equation}\label{eq: D}
D=\bigcup_{j\in\N}\(g_j+\vep  N_j^{\circ}\)\supset C \supset F.
\end{equation}
We show first that $D$ is a $\omega^*$-compact set. Let $\{h_m\}_{m\in\N} \subset D$,
$h_m \overset {\omega^*}\to h_0\in 2B_{X^*}$. If infinitely many of $h_m$ are in
one of the sets $g_j+\vep N_{j}^{\circ}$, then of course $h_0 \in D$. So let us
assume that $h_m \in g_{j_m}+ \vep N_{j_m}^{\circ}$ where $m\in\N$,
$j_m\rightarrow \infty$ as $m\rightarrow \infty$. Then, by \eqref{eq: (c)} we may write $h_m=g_{j_m}+z_{m}+\vep u_m$
where $u_m\in B_{M_{n_{j_{m}}}}$, $m\in\N$, and $\|z_{m}\|\to 0$ as $m\to \infty$.
Since $j_m\rightarrow \infty$, we have $n_{j_m} \rightarrow \infty $, according to Lemma \ref{l: decomp}, and therefore
$\omega^* -\lim u_m= 0$. Thus $h_0=\omega^*-\lim g_{j_m} \in F \subset D$,
which completes the proof of $\omega^*$-closedness of the set $D$.

Write $Q^\circ=\overline{\conv}^{\omega^{\ast}}(D)$ and $Q=\{x\in X:\ f(x)\leq 1\ \mathrm{for\ all}\ f\in Q^\circ\}$.
The latter set is a closed convex set in $X$ with a non-empty interior and it is used here to approximate $W$.
Let us denote by $F_Q$ the Minkowski functional of the set $Q$.

For each $j\in\N$ we define
\begin{equation}\label{eq: fksup}
F_j(x)=\sup\{f(x):\ f\in g_j+\vep N_{j}^{\circ}\}=g_j(x)+\vep N_{j}(x).
\end{equation}

Observe that $F_j$ are positively homogeneous, $2$-Lipschitz for  small enough $\vep$,
and by the right hand side they are $C^k$-smooth away from the origin.

Pick $x\in Q$ such that $\sup_{f\in Q^\circ}f(x)=1$. Then necessarily $\sup_{f\in D}f(x)=1$.
By the $\omega^*$-compactness of $D$ there is $f\in D$ such that $f(x)=1$ and of course $f(\lambda x)=\lambda$
for $\lambda \in \R$. This means, by recalling \eqref{eq: D} and \eqref{eq: fksup}, that
\[F_Q=\max_j F_j.\]

We introduce a strongly lattice function $F\colon \ell_\infty\to [0,\infty]$ as follows
\begin{equation}\label{eq: F}
F(x_k)=\sum\phi_k(x_k),
\end{equation}
where $\phi_j:\R\to[0,\infty)$ are $C^\infty$-smooth convex and even functions such that

\begin{equation}\label{eq: phidef}
\phi_k[-1-\frac{k}{k+2}\delta_k, 1+\frac{k}{k+2}\delta_k]=0,\ \phi_k(1+\frac{k+1}{k+2}\delta_k)=1, \phi_k(1+\delta_k)=3.
\end{equation}

Let $0<\lambda_{1}<1$. Fix a decreasing sequence $\delta_n\to 0$, such that
\begin{equation}\label{eq: delta}
\lambda_{n}=\frac{1+\frac{n}{n+2}\delta_{n}}{1+\delta_{n}}=1-\frac{2\delta_{n}}{(n+2)(1+\delta_{n})}
\end{equation}
is an increasing sequence and equals to $\lambda_{1}$ for $n=1$.
Let $h_k=(1+\delta_k)F_k$.
By Lemma \ref{l:SLconvexComposition} we obtain that

$$
G(x)=\sum\phi_k(h_k(x))
$$
is convex wherever it is well-defined. Next, we will check that the formula for $G$ locally contains only finitely many non-zero summands on the set $G^{-1}[0,1]$. Indeed, let $x$ be such that $F_Q(x)=1$, and suppose $F_k(x)=1$.
If $G(\lambda x)=1$ then $\phi_k\circ h_k(\lambda x)\le 1$, thus $h_k(\lambda x)\le 1+\frac{k+1}{k+2}\delta_k$, so
$F_k(\lambda x)\le \frac{1+\frac{k+1}{k+2}\delta_k}{1+\delta_k}$, and hence
\begin{equation}\label{eq: lambdafrac}
\lambda\le \frac{1+\frac{k+1}{k+2}\delta_k}{1+\delta_k}.
\end{equation}

Put $\rho=\delta_{k}/(4(k+2)(1+\delta_{k}))$. Now, if $\|y-\lambda x\|<\rho$, then by using the facts that
$F_{n}$ is $2$-Lipschitz, $F_{Q}(x)=1$ and \eqref{eq: lambdafrac} we obtain

\[F_n(y)\le \lambda F_n(x) +2\rho\le \lambda+2\rho\le \frac{1+\frac{k+1}{k+2}\delta_k}{1+\delta_k}+2\rho=1-2\rho\quad \mathrm{for}\ n.\]

Thus
\[h_N(y)\le(1+\delta_N)(1-2\rho)\le 1+\frac{N}{N+2}\delta_N,\]
provided that $N$ is large enough. By recalling \eqref{eq: phidef} we observe that $G(y)$ is a finite sum for
$\|y-\lambda x\|<\rho$. Thus $G$ is $C^k$-smooth.

Finally, let us check that $G^{-1}([0,1])$ is close to $Q$. Indeed, if $|F_{Q}(x)|=1$, then $G(x)\geq 3$ according to
\eqref{eq: phidef}, so that $G^{-1}([0,1])\subset Q$. If $x\in Q$ is such that $|F_{n}(x)|\leq \lambda_{n}$ for $n$,
then $G(x)=0$. Since $\lambda_{n}$ is an increasing sequence, we get that
\[\{x\in Q:\ F_{Q}(x)\leq \lambda_{1}\}\subset G^{-1}(0).\]
Since $\lambda_{1}<1$ was arbitrary, the proof is complete.
\end{proof}

 Similarly to \cite{defoha-separ} we conclude the following fact.

\begin{corollary}
Let $(X,\|\cdot\|)$ be a separable normed space that satisfies the
assumptions of any of the following theorems contained in the above mentioned paper:
Corollary 2.4, Corollary 3.10, Theorem 4.1 and Theorem 4.3.

Then arbitrary convex function on $X$ can be approximated on bounded sets by
convex functions of the same degree of smoothness as are the approximations
of norms in the corresponding theorem.
\end{corollary}


\begin{thebibliography}{DGZ}

\bibitem[Ben]{ben} Y. Benyamini, {\em An extension theorem for separable Banach
spaces,} Israel J. Math., Vol. 29, No. 1, 1978, 24-30.\par


%\bibitem[BCR]{bcr} J. Bochnak, M. Coste, M. Roy, {\em Real Algebraic Geometry}, Springer, 1998.\par


\bibitem[BoFr]{bonfram} R. Bonic and J. Frampton, {\em Smooth functions on Banach manifolds}, J. Math.\ Mech.~{\bf 15} (1966), 877--898.


\bibitem[CH]{ch} M. Cepedello-Boiso and P. H\' ajek, {\em Analytic approximations of uniformly
continuous functions in real Banach spaces},  J. Math. Anal. Appl. 256
(2001), 80--98.



\bibitem[BoWa]{borvein-vander} J. Borwein and J. Vanderwerff,
{\em Convex Functions: Constructions, Characterizations and Counterexamples}.
Encyclopedia Math. Appl. vol. 109, Cambridge University Press, 2009.


\bibitem[D]{D} R. Deville, {\it Geometrical implications of the existence of very
smooth bump functions in Banach spaces,} Israel J. Math., 6(1989), 1-22.\par

\bibitem[DFH1]{defoha-separ} R. Deville, V.P. Fonf and P. H\'ajek, {\it
Analytic and $C^k$ approximations of norms in separable Banach spaces}, Studia
Math.~{\bf 120} (1996), 61--74.

\bibitem[DFH2]{defoha-polyhedral} R. Deville, V.P. Fonf and P. H\'ajek, {\it
Analytic and polyhedral approximation of convex bodies in separable polyhedral
Banach spaces}, Israel J. Math.~{\bf 105} (1998), 139--154.


\bibitem[DGZ]{dgz} R. Deville, G. Godefroy, and V. Zizler, {\it Smoothness and
renormings in Banach spaces,} Pitman Monographs and Surveys in Pure and
Applied Mathematics, 64, 1993.\par


\bibitem[Dieu]{dieudonne} J. Dieudonn\' e, {\em Foundations of Modern Analysis}, Academic Press, 1960.

\bibitem[FHHMZ]{fhhmz} M. Fabian, P. Habala, P. H\'ajek, V. Montesinos, V. Zizler,
{\it Banach Space Theory, The Basis for Linear and Nonlinear Analysis}, CMS Books in Mathematics, Springer 2011.\par

\bibitem[FPWZ]{fpwz} M. Fabian, D. Preiss, J.H.M. Whitfield and V. Zizler, {\it
Separating Polynomials on Banach Spaces,} Quart. J. Math. Oxford (2), 40,
1989, 409--422.\par

\bibitem[Fe]{fe} H. Federer, {\it Geometric theory of measure,} Springer-Verlag,
1970.\par

\bibitem[Fry1]{fry1} R. Fry, {\em Analytic approximation on $c_0$}, J. Funct. Anal. 158 (1998), 509--520.

\bibitem[Haje1]{hajek-locally} P. H\'ajek, {\em Smooth norms that depend
locally on finitely many coordinates}, Proc.\ Amer.\ Math.\ Soc.~{\bf 123}
(1995), 3817--3821.

\bibitem[HMVZ]{bos} P. H\'ajek, V. Montesinos, J. Vanderwerff, and V. Zizler, {\em
Biorthogonal systems in Banach spaces}, CMS Books in Mathematics, Canadian
Mathematical Society, Springer Verlag, 2007.


\bibitem[HP]{hp} P. H\' ajek and A. Proch\' azka,
{\em $C^k$-smooth approximations of LUR norms}, to appear in Trans. AMS.


\bibitem[Kurz]{kurzweil} J. Kurzweil, {\em On approximation in real Banach
spaces}, Studia Math.~{\bf 14} (1954), 213--231.


%\bibitem[Sz]{szpond} J. Szpond,
%{\em Reduction of power series in a polydisc with respect to a Gröbner basis},
%Bull. Pol. Acad. Sci. Math. {\bf 53} (2005), 137-145.\par

\bibitem[Zp]{zp} M. Zippin, {\it The separable extension problem,} Isr. J. Math.,
Vol. 26, 1977, 372-387.\par


\end{thebibliography}
\end{document}